# A New Adaptive Channel Estimation for Frequency Selective Time Varying Fading OFDM Channels

Wessam M. Afifi, Hassan M. Elkamchouchi

*Abstract—* In this paper a new algorithm for adaptive dynamic channel estimation for frequency selective time varying fading OFDM channels is proposed. The new algorithm adopts a new strategy that successfully increases OFDM symbol rate. Instead of using a fixed training pilot sequence, the proposed algorithm uses a logic controller to choose among several available training patterns. The controller choice is based on the cross-correlation between pilot symbols over two consecutive time instants (which is considered to be a suitable measure of channel stationarity) as well as the deviation from the desired BER. Simulation results of the system performance confirm the effectiveness of this new channel estimation technique over traditional non-adaptive estimation methods in increasing the data rate of OFDM symbols while maintaining the same probability of error.

## I. INTRODUCTION

Orthogonal Frequency Division Multiplexing (OFDM) is a multicarrier modulation scheme which is based on the idea of dividing a given high-bit-rate serial data stream into several parallel lower bit-rate streams and modulating each stream on separate orthogonal carriers, often called subcarriers [1]. Thereby, the bandwidth of the subcarriers becomes small compared with the coherence bandwidth of the channel; that is, the individual subcarriers experience flat fading, which allows for simple equalization. So, OFDM solves the problem of the frequency selective channels as it divides the total bandwidth into number of subcarriers each of them have small bandwidth and will see the channel as flat. This implies that the symbol period of the substreams is made long compared to the delay spread of the time-dispersive radio channel as the symbol duration is inversely proportional to the data rate. This will help to eliminate or minimize the Inter Symbol Interference (ISI).

OFDM has a lot of applications and is used in many systems as (ADSL) Asymmetric Digital Subscriber Line[5], (PLC) Power Line Communication, IEEE802.11 (WLAN), IEEE802.16 (WiMAX), (DAB) Digital Audio Broadcasting systems and (DVB) Digital Video Broadcasting [4]

In this paper we propose a new algorithm for OFDM time varying channel estimation. The algorithm further enhances the performance of OFDM-based systems by successfully combining the advantages of high symbol rate and low BER through the adoption of an adaptive scheme. The proposed adaptive scheme modifies the sent training pattern of pilot symbols according to the channel status.

The rest of the paper is organized as follows. In section II, we give a brief overview of the basic concepts behind the different channel estimation approaches. This will help us compare and contrast their advantages and disadvantages in the following section. In section III, we discuss the new proposed adaptive channel estimation technique. The technique is illustrated by discussing the adopted system model and algorithm flow chart in details. Section IV shows the results of different comparative simulation case studies. The results clearly highlight the advantages of the proposed technique. Finally, section V concludes the paper by summarizing our contributions and briefly discussing directions for future research.

## II. CHANNEL ESTIMATION APPROACHES

Channel estimation has a great importance in any wireless communication system. When the information is sent from the transmitter to the receiver, it suffers from amplitude scaling and phase rotation. As the wireless communication channel may be frequency selective time varying fading channel so the data absolutely will suffer from great distortion. Time varying fading channels has an impulse response which is a function of the time, i.e. the channel impulse response differs with the time. While frequency selective fading channels have an impulse response which differ from a band of frequency to another or in other words every frequency or band of frequencies suffers from different fading. Thus, channel estimation is done in the wireless communication systems to compensate for phase rotation and amplitude scaling. Channel estimation in wired communication systems differ from the wireless one as the channel is estimated at start-up, and since channel remains the same, therefore no need to estimate it continuously [2].

Channel estimation has mainly two types which are Pilot Symbol Assisted Modulation (PSAM) and Blind channel estimation. Channel estimates in the PSAM are often achieved by multiplexing known symbols, so called, pilot symbols into the data sequence (useful information) for the purpose of channel sounding. Pilot symbols are transmitted at certain locations of the OFDM frequency time lattice, instead of data. General fading channel can be viewed as a 2-D signal (time and frequency), which is sampled at pilot positions and channel attenuations between pilots are estimated by interpolations [3].

The problem with this type of channel estimation is to decide where and how often to insert pilot symbols. The number of pilot tones necessary to sample the transfer function can be determined on the basis of sampling theorem as follows. The frequency domain channel's transfer function $H(f)$ is the fourier transform of the impulse response $h(t)$ [2]. Each of the impulses in the impulse response will result a complex exponential function $e^{\frac{j2\Pi\tau}{T_s}}$ in the frequency domain, depending on its time delay $\tau$, where $T_s$ is the symbol time. In order to sample this contribution to $H(f)$ according to the sampling theorem, the maximum pilot spacing $\Delta p$ in the OFDM symbol is [2]:

$$\Delta p \leq \frac{N}{2\tau/T_s} \Delta f \quad (1)$$

Where $\Delta f$ is the subcarrier bandwidth.

PSAM can be achieved by either inserting pilot tones into all of the subcarriers of OFDM symbols with a specific period (block type pilot channel estimation) or inserting pilot tones into each OFDM symbol (comb type pilot arrangement).
In fig. 1, the block type pilot channel estimation has been developed under the assumption of slow fading channel or when the channel transfer function is not changing very rapidly. The comb type pilot channel estimation has been developed for fast varying channels as demonstrated in Fig. 2. The second type for channel estimation is the blind channel estimation where the receiver must determine the channel without the aid of known symbols. Blind means that pilot symbols are not available to the receiver [1]. This type uses some information about the channel statistics to estimate the channel.

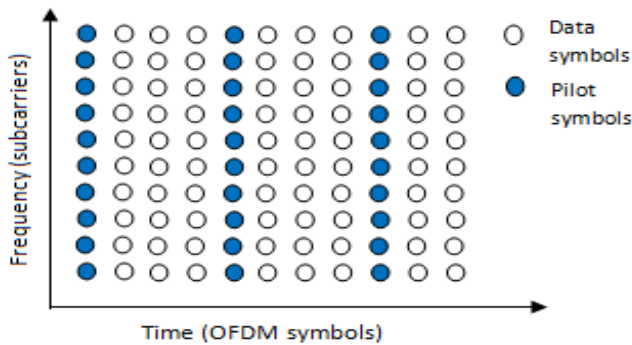

Fig. 1. Block Pilot Patterns

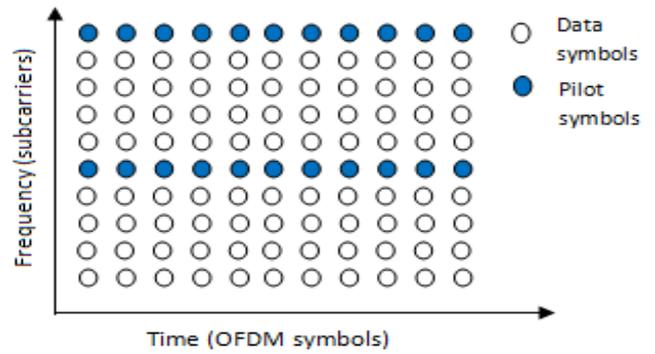

Fig. 2. Comb Pilot Patterns

III. PROPOSED APPROACH

It is clear that sending some pilot symbols instead of the OFDM symbols in the first type of the channel estimation (PSAM) for the purpose of channel estimation, results in a decrease of OFDM symbols rate. This decrease in the data rate of the OFDM symbols is accompanied, however, by a decrease in probability of error. This is due to the high accuracy of the obtained channel estimate.

On the other hand, in the second type of the channel estimation (blind channel estimation), no pilot symbols are sent to estimate the channel. Instead, statistical information about the channel is used for the purpose of channel estimation, which results in higher data rate compared with PSAM. However, the probability of error with blind channel estimation techniques is higher than that with PSAM. Moreover, blind techniques add computational complexities and may be difficult to implement in case of real time applications.

The new adaptive channel estimation technique proposed in this paper aims to combine the advantages of both approaches of channel estimation.. The objective of this new algorithm is to increase the data rate of the OFDM symbols for the same probability of error for any type of PSAM.

Comparing this new algorithm with blind channel estimation will clearly show that the proposed algorithm gives lower probability of error at the expense of a slightly lower data rate. The proposed algorithm is best described by the system model depicted in Fig. 3. We start the communication process by sending the preamble and the pilot symbols for the purpose of channel sounding multiplexed with the information data.
Afterwards, the cross-correlation between pilot symbols at two consecutive time instants is calculated. First, we multiply the two pilot symbols $P_1(t)$ and $P_2(t)$, then integrate the output so that we can get the cross correlation between these two pilot symbols. The cross correlation $R_{P1P2}$ can be calculated from equation (2):

$$R_{P1P2} = \frac{1}{\sqrt{E_1 E_2}} \int_0^T P_1(t) P_2(t) dt \quad (2)$$

Where $P_1(t)$ and $P_2(t)$ are the two pilot symbols.

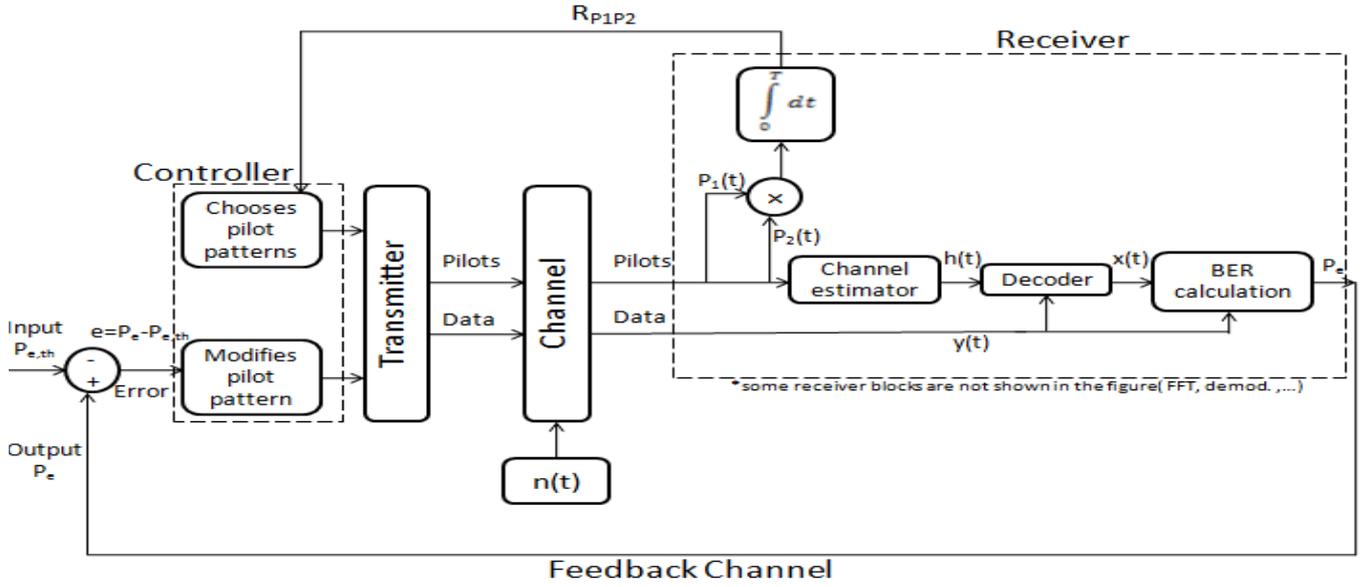

Fig. 3. System model with the communication process and controller

The value of the cross correlation is then sent back to the controller at the transmitter. The controller chooses the corresponding pattern of pilot symbols according to the value of the cross correlation. At the receiver, the pilot symbols are used for the purpose of channel estimation, so that we can get the channel impulse response $h(t)$. Using the received data and the estimated channel impulse response we can reconstruct the transmitted signal and calculate the probability of error (BER). The probability of error is then sent back to the transmitter through the feedback channel. At the transmitter, we compare between the threshold or reference probability of error ($P_{e,th}$) of error and that calculated ($P_e$).

$$e = P_e - P_{e,th} \tag{3}$$

This error $e$ is calculated at the instants of not sending pilot symbols. The pilot pattern chosen based on the cross-correlation value is maintained until this error $e$ becomes positive. In this case, the "modify pattern" block of the controller is activated. This block re-triggers the procedure of channel sounding/cross-correlation calculations and modifies the current pilot pattern accordingly. We notice in the system model that the rate of the feedback channel is low. We use the feedback channel either to send the value of the cross correlation or to send the value of the probability of error. Thus, (as will be further clarified through results), the required feedback rate does not impair the operation of the system or adversely affect its superior performance.

The details of the strategy employed by the controller to choose a suitable pilot pattern are explained in the flow chart in Fig. 4. The controller checks the value of the cross correlation between pilot symbols and accordingly estimates the stationarity of the channel. According to the value of the cross correlation, the algorithm chooses the pattern of the pilot symbols. Many simulations have been performed to find the optimal boundaries for the values of the cross correlation, i.e. for each boundary, we choose a specific pattern of pilot symbols. Table I summarizes the decided boundaries/patterns configurations based on these simulation results. Based on Table I, the controller has been designed and a general formula for defining the boundaries of the cross-correlation for each pattern has been developed. The lower boundaries (LB) of the cross correlation values are 0, 0.7, 0.8 and 0.9. They can be described by equation (4).

$$LB = 0.7i - 0.3i(i-1) + 0.1i(i-1)(i-2) \tag{4}$$

Where (i) takes values from 0 to 3 which are the numbers of the boundaries (4 boundaries) $0 \le i \le 3$

$$LB = 0.1i^3 - 0.6i^2 + 1.2i \tag{5}$$

The higher boundaries (HB) of the cross correlation values are 0.7, 0.8, 0.9 and 1. These are described by (6).

$$HB = 0.7(i+1) - 0.6i \tag{6}$$

$$HB = 0.1i + 0.7 \tag{7}$$

$$LB \le R_{P1P2} < HB \tag{8}$$

$$0.1i^3 - 0.6i^2 + 1.2i \le R_{P1P2} < 0.1i + 0.7 \tag{9}$$

TABLE I
CROSS CORRELATION VALUES WITH THE APPROPRIATE PATTERNS OF PILOT SYMBOLS

| Cross-correlation boundaries | Pattern used | No. of pilot symbols sent |
|---|---|---|
| [0,0.7[ | Pattern (1) | L pilot symbols every n time instances |
| [0.7,0.8[ | Pattern (2) | L pilot symbols every 2n time instances |
| [0.8,0.9[ | Pattern (3) | L pilot symbols every 4n time instances |
| [0.9,1] | Pattern (4) | L pilot symbols every 8n time instances |

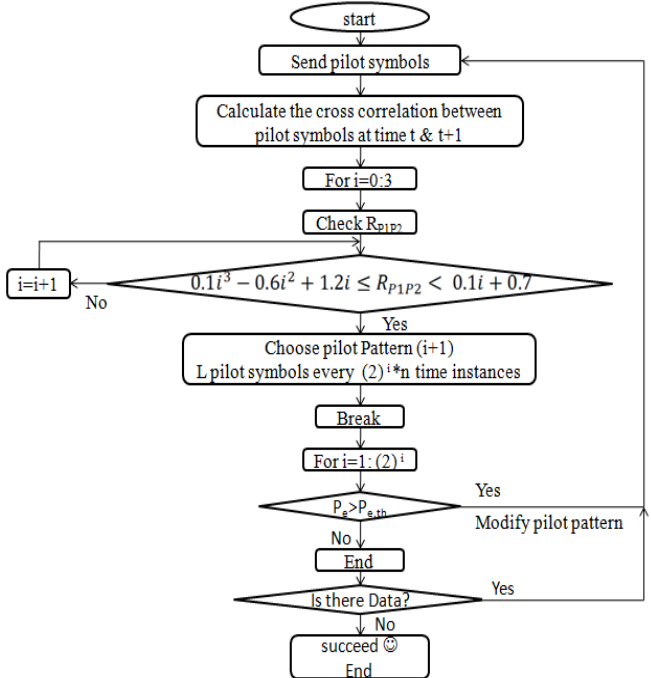

Fig. 4. Flowchart of the proposed algorithm

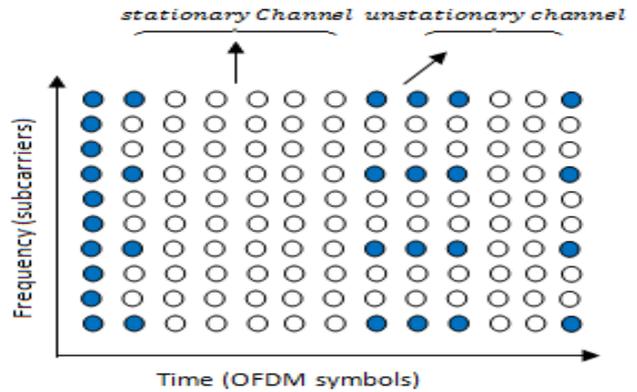

Fig. 5. Adaptive pilot pattern (Example 1)

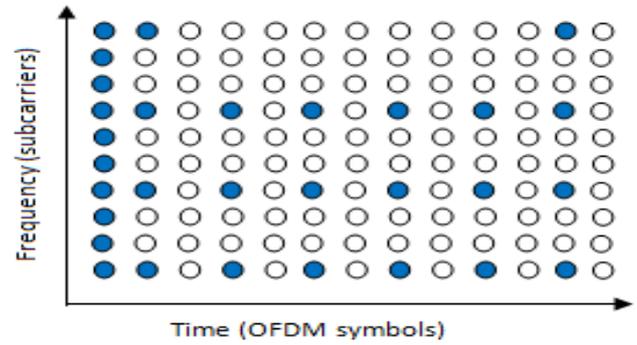

Fig. 6. Adaptive pilot pattern (Example 2)

According to the value of the cross correlation and the index ($i$), the controller iterates based on (9) to determine the used pattern of pilot symbols. The controller sends pilot symbols every $2^i \times n$ time instances, where ($n$) is a parameter that depends on the stationarity of the channel. This parameter, which may be set to one, can be varied to tune the performance of the algorithm. During the period of sending data symbols, we calculate the probability of error. If the calculated probability of error is larger than the threshold probability of error, the controller resets the algorithm back to pattern (1).

Figs. 5 and 6 show how the proposed algorithm senses the channel nature and adapts the pilot patterns accordingly. Fig. 5 shows an example of the adaptive channel estimation where the algorithm exploits the nature of the channel. When the channel is stationary the algorithm plays a role in not sending pilot symbols and send OFDM (useful information) instead and the receiver could use the same previous pilot symbols for estimating the channel as the channel seems to be stationary.

When the channel becomes unstationary or fast varying the algorithm increases the number of pilot symbols for the purpose of channel sounding or for intensive channel estimation.

As a result this new algorithm tracks the channel and chooses the appropriate pattern of pilot symbols depending on the characteristics of the channel which helps in saving the pilot symbols and increasing the data rate of OFDM symbols.

Another Example for the Adaptive channel estimation is shown in Fig. 6. where the algorithm calculates the cross correlation between pilot symbols and recognizes that the higher band of frequency of the channel is stationary. The algorithm decided to decrease the number of pilot symbols and use the previous pilot symbols in estimating the channel in this higher band only. For the other bands of frequency the channel varies rapidly with time so we increase the number of pilot symbols to estimate the channel correctly. This example, however, is not currently implemented but is included to show an exciting future direction for the proposed approach.

IV. SIMULATION RESULTS

Simulation case studies have been performed to evaluate the proposed algorithm performance. In all case studies the new adaptive algorithm is applied to the least square channel estimation and compared with it. Notice that the new adaptive channel estimation can be applied to any type of channel estimation as minimum mean square error.

Three types of simulations have been performed:

In the first type of simulations, we show the effectiveness of the new proposed algorithm over traditional channel estimation. In the second: for a specific model of the channel, we perform simulations for different values of the LB and the HB to find the optimal boundaries for cross-correlation/ patterns configuration.

In the Third set of simulations performed, the LB and HB values are kept constant at their optimal boundaries determined through the second set of simulations. The aim of these simulations is to evaluate the performance of the proposed algorithm for different simulated channel models.

Results of the first part of the simulation case studies are shown in Figs. 7 and 8. Fig. 7 shows the relation between the BER and the SNR in dB. It shows that the probability of error of the adaptive channel estimation "the new algorithm" is nearly the same as the probability of error of the LS channel estimation "the traditional estimation".

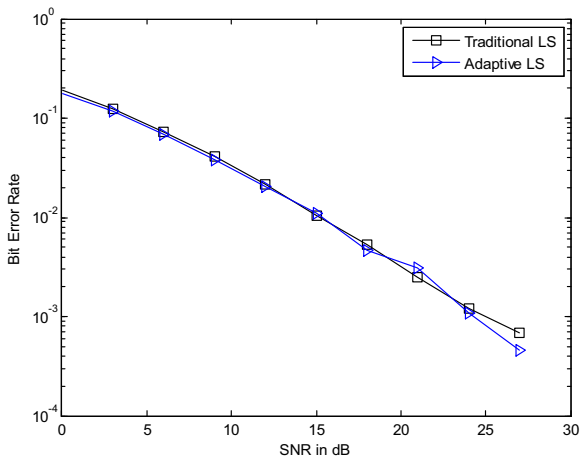

Fig. 7. BER vs SNR in dB

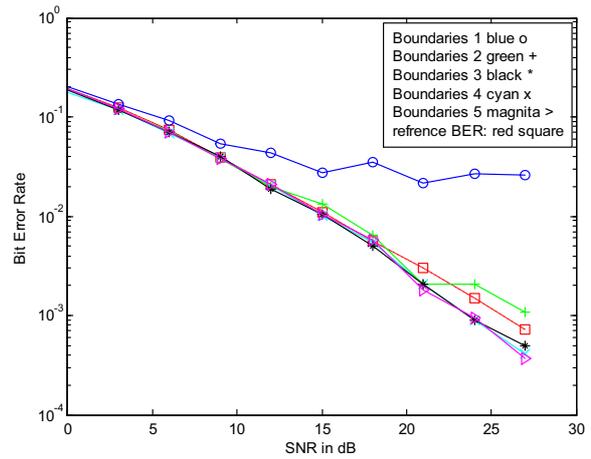

Fig. 9. BER vs SNR in dB

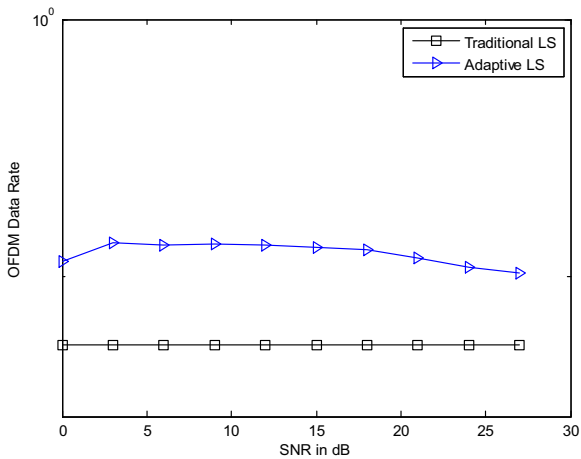

Fig. 8. OFDM data rate vs SNR in dB

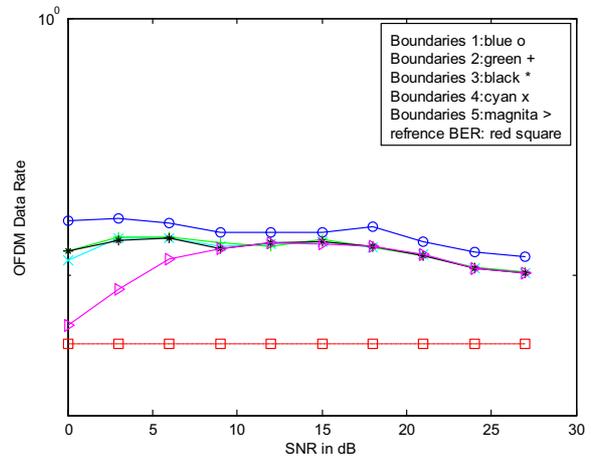

Fig. 10. OFDM data rate vs SNR in dB

Fig. 8 shows the relation between the data rate of the OFDM systems and the SNR in dB. It shows that by using the adaptive channel estimation we get a higher data rate than that of traditional channel estimation. The reason for this increase in the data rate is that at certain time instances, we send OFDM data symbols instead of sending pilot symbols. Figs. 7 and 8 clearly show the ability of the adaptive channel estimation to track the OFDM channel. Moreover, the algorithm maintains nearly the same threshold probability of error. When the BER increases above the threshold BER, the algorithm increases the number of pilot symbols to compensate for this increase and hence the data rate decreases slightly.

The second part of the simulation case studies show the parametric study that we performed for different values of the lower boundary (LB) and the higher boundary (HB) for a specific model of the channel. The aim of this set of simulations is to determine the optimal values for the LB and HB The value of the cross-correlation here can take values from 0 to 1, so we simulate for some boundaries of the cross-correlation. The results are shown in Figs. 9 and 10.

TABLE II
CROSS CORRELATION VALUES FOR DIFFERENT BOUNDARIES

| Boundary pattern | 1 | 2 | 3 | 4 | 5 |
|---|---|---|---|---|---|
| Pattern1 | [0,0.25[ | [0,0.5[ | [0,0.6[ | [0,0.7[ | [0,0.9[ |
| Pattern2 | [0.25,0.5[ | [0.5,0.7[ | [0.6,0.7[ | [0.7,0.8[ | [0.9,0.95[ |
| Pattern3 | [0.5,0.75[ | [0.7,0.9[ | [0.7,0.9[ | [0.8,0.9[ | [0.95,0.97[ |
| Pattern4 | [0.75,1] | [0.9,1] | [0.9,1] | [0.9,1] | [0.975,1] |

Table II shows the values of the cross-correlation for different boundaries. Fig. 9 shows the BER of the different boundaries. It is shown from the figure that although boundary (1) has higher data rate it is refused as it has high BER. Boundary (2) is refused also as it gives high BER. Boundary (5) gives a lower data rate than other boundaries so it is refused. Boundary (3) and (4) nearly the same and can give almost the same results however boundary (4) has a slightly lower BER than boundary (3). Boundary (4) is the optimal one which gives the optimal values for the LB and HB.

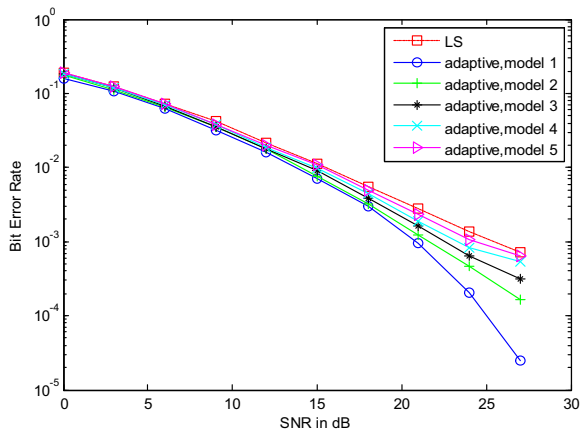

Fig. 11. BER vs SNR in dB

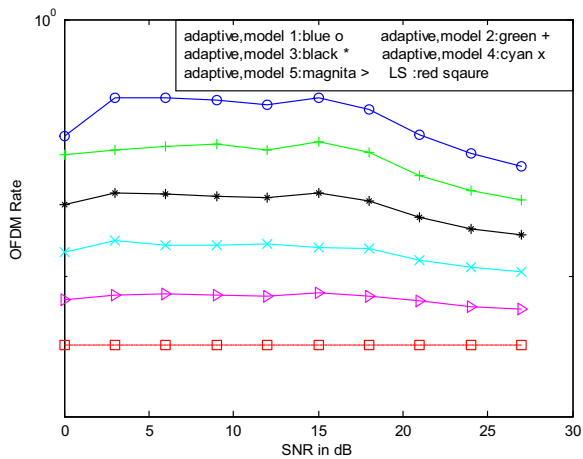

Fig. 12. OFDM data rate vs SNR in dB

In the third set of simulations we use the optimal values for the LB and HB (Boundary 4). The aim of this set of simulations is to evaluate the performance of the proposed algorithm for different simulated channel models.

TABLE III
DESCRIPTION OF CHANNEL MODELS USED

| Model Description | 1 | 2 | 3 | 4 | 5 |
|---|---|---|---|---|---|
| Percentage of channel stationarity | 100% | 80% | 60% | 40% | 20% |
| Percentage of channel randomization | 0% | 20% | 40% | 60% | 80% |

Table III shows the description of the different models used in this set of simulations.

Fig. 11 shows that all different channel models nearly give the same BER however the algorithm works better for channels with high percentage of stationarity. Notice that all curves for the BER are within the range of the threshold BER.

In addition to the aforementioned simulation case studies, Fig. 12 is included to add further insight to another aspect of the proposed algorithm performance. The figure clearly depicts that the algorithm exhibits "graceful degradation". This means although the algorithm performs better for channels with a higher degree of stationarity, its performance at lower degrees of stationarity is still superior (in a multi-disciplinary sense) to other channel estimation alternative approaches.

## V. CONCLUSION

The main contributions of the paper is offering a new innovative adaptive scheme for pilot patterns selection together with defining new simple channel state measures (i.e. channel degree of stationarity is estimated based on cross-correlation, while the algorithm performance as a whole is estimated based on deviation from desired BER) to implement such a scheme using a suitable controller. Simulation case studies clearly demonstrated that the proposed technique simultaneously combines the advantages of low BER and high data rate.

Several interesting directions for future work exist. First, it is possible to employ the algorithm both in time and frequency domains. Second, the crisp logic controller may be replaced by a fuzzy logic controller in which case a tradeoff between reliability and transmitter complexity should be evaluated. Furthermore, instead of using the deviation from desired BER as a criterion for deciding whether or not to reset the algorithm, a parametric study may be carried out (similar to that performed to set cross-correlation optimal boundaries) to determine optimal error boundaries for different groups of patterns.

ACKNOWLEDGMENT

Thanks to Allah almighty whose love and desire to satisfy motivated this work and every good we do. Thanks to dr. Ahmed K. Sultan for his great effort and contribution in this work

REFERENCES

[1] J. Andrews, A. Ghosh, R. Muhamed, *Fundamentals of WiMAX,* Prentice Hall Communications Engineering and Emerging Technologies Series, 2007.
[2] K. Arshad "Channel Estimation in OFDM Systems", 6th Saudi Engineering Conference, Dhahran, Saudi Arabia, Dec. 2002 .
[3] F. Tufvesson, T. Maseng, "Pilot Assisted Channel Estimation for OFDM in Mobile Cellular Systems," Proceedings of IEEE Vehicular Tech. Conference, Phoenix USA, pp. 1639–1643, May 1997.
[4] Henrik Schulze and Christian Theory and Applications of OFDM and CDMA: Wideband Wireless Communications, Wiley, NY, U.S.A, 2005
[5] I. Bradaric and A.P. Petropulu, "Blind Estimation of the Carrier Frequency Offset in OFDM Systems," *4th IEEE Workshop on Signal Processing Advances in Wireless Communications, SPAWC'03*, Rome, Italy, June 2003.
[6] T.Wang;,J.Proakis,J Zeidler, "Techniques for Suppression of Intercarrier Interference in OFDM systems", *IEEE Wireless Communications and Networking Conference*, Page(s): 39 - 44 Vol. 1, March 2005
[7] L.J. Cimini, Jr., "Analysis and Simulation of a Digital Mobile Channel using Orthogonal Frequency Division Multiplexing," in IEEE Trans. Commun. Technol., vol. 33, pp. 665-675, July 1985.